\documentclass[a4paper,12pt]{amsart}
\usepackage{fullpage, verbatim}
\usepackage{amsfonts, amssymb, amsthm, mathrsfs}
\usepackage[shortlabels]{enumitem} 
\usepackage{tikz}
\usepackage{tikz-cd}
\usetikzlibrary{shapes,decorations.pathmorphing,calc,arrows}

\usepackage[style=alphabetic]{biblatex}
\usepackage{caption} 

\tikzset{%
  dots/.style args={#1per #2}{%
    line cap=round, dash pattern=on 0 off #2/#1 },
  DyckWidth/.style={rounded corners=1, line width=1.5pt}
}

\definecolor{lightgrey}{rgb}{0.7,0.7,0.7}

\makeatletter
\long\def\ifnodedefined#1#2#3{%
    \@ifundefined{pgf@sh@ns@#1}{#3}{#2}%
}
\makeatother

\newlength{\unit}
\usepackage[colorinlistoftodos]{todonotes}
\PassOptionsToPackage{pdfusetitle,colorlinks}{hyperref}
\usepackage{bookmark}
\hypersetup{
  linkcolor={red!70!black},
  citecolor={green!70!black},
  urlcolor={teal!90!black}
}
\numberwithin{equation}{section}
\usepackage{cleveref}

\newcommand{\Dfn}[1]{\emph{\color{blue}#1}} 

\numberwithin{equation}{section}
\theoremstyle{plain}
\newtheorem{lemma}[equation]{Lemma}
\newtheorem{problem}[equation]{Problem}
\newtheorem{theorem}[equation]{Theorem}

\newtheorem{corollary}[equation]{Corollary}
\newtheorem{proposition}[equation]{Proposition}

\theoremstyle{definition}
\newtheorem{definition}[equation]{Definition}
\newtheorem{remark}[equation]{Remark}

\newtheorem{example}[equation]{Example}

\DeclareMathOperator{\pdim}{pdim}%
\DeclareMathOperator{\idim}{idim}%
\DeclareMathOperator{\gldim}{gldim}%
\newcommand{\drawPath}[4]
{%
  \draw[rounded corners=1, line width=2, #3] #4
  \foreach \dir in {#1}%
  {
    \ifnum\dir=#2
    -- ++(1,0)
    \else
    -- ++(0,1)
    \fi
  };}

\title[Nakayama algebras of small homological dimension and pattern avoiding permutations]%
{Nakayama algebras of small homological dimension and pattern avoiding permutations}%

\date{\today}
\author[V.Kl\'asz]{Vikt\'oria Kl\'asz$^\dagger$}%
\address[V.~Kl\'asz]{Mathematical Institute of the University of Bonn, Endenicher Allee 60, 53115 Bonn, Germany}%
\email{klasz@math.uni-bonn.de}%
\thanks{$^\dagger$Supported by the Deutsche Forschungsgemeinschaft
(DFG, German Research Foundation) under Germany’s Excellence Strategy - GZ 2047/1, Projekt-ID
390685813}

\author[R.~Marczinzik]{René Marczinzik}%
\address[R.~Marczinzik]{Mathematical Institute of the University of Bonn, Endenicher Allee 60, 53115 Bonn, Germany}
\email{marczire@math.uni-bonn.de}

\author[J. Marquardt]{Judith Marquardt}
\address[J.~Marquardt]{Mathematical Institute of the University of Bonn, Endenicher Allee 60, 53115 Bonn, Germany}
\email{judith.marquardt@uni-bonn.de}
\subjclass[2010]{Primary 16G10, 16E10,05E10}
\keywords{Nakayama algebras, Dyck paths, shod algebras, pattern avoiding permutations}

\begin{document}
\begin{abstract}
We give a combinatorial classification of Nakayama algebras of small homological dimension using the Krattenthaler bijection between Dyck paths and 132-avoiding permutations.
\end{abstract}

\maketitle
\setcounter{tocdepth}{1}

\section{Introcution}
Following the work of \cite{CL} and \cite{RS}, we call a finite-dimensional algebra $A$ \Dfn{shod} (for \textbf{s}mall \textbf{ho}mological \textbf{d}imension) if the projective or injective dimension of every indecomposable $A$-module $X$ is at most 1.
Recently, \cite{BZ} showed that shod algebras can be characterised using silting theory. 
It is known that such algebras have global dimension at most 3. Shod algebras of global dimension 3 are called \Dfn{strict shod}. Those with global dimension at most 2 coincide with \Dfn{quasi-tilted algebras}, see \cite[Chapter II]{HRS} for more information.
Quasi-tilted algebras contain the class of \Dfn{tilted algebras} defined as the endomorphism algebras of tilting modules over hereditary algebras which play a central role in the representation theory of finite-dimensional algebras, see for example the textbooks \cite{ASS,SS1,SS2}. Tilted algebras and their relation extension algebras, called cluster tilted algebras, also play an important role in the theory of cluster algebras, see for example \cite{A,BMRRT}.
A finite-dimensional algebra $A$ is called a \Dfn{Nakayama algebra} if every indecomposable $A$-module is uniserial, meaning it has a unique composition series. The representation theory of Nakayama algebras is a classical topic, see for example the textbooks \cite{AF,ASS,SY,ARS}. We remark that they are also often called serial algebras in the literature, such as in the textbook \cite{AF}.
Throughout this article, we will assume that every Nakayama algebra is given by a quiver and relations. Then the simple modules correspond to the vertices of the quiver. If we work over an algebraically closed field, this allows us to obtain, up to Morita equivalence every Nakayama algebra, see for example \cite[Chapter II, Theorem 3.7]{ASS}.
We can also assume that the algebras have a connected quiver and so unless otherwise stated we will assume that Nakayama algebras are connected.
Nakayama algebras can be grouped into two classes: Those with a linear quiver and those with a cyclic quiver. It is known, see \cite{CL}, that shod algebras always have an acyclic quiver and therefore, we will always work with linear Nakayama algebras in this article.
Recently, it was noted that linear Nakayama algebras with a fixed number of simple modules correspond in a canonical way to Dyck paths via their Auslander-Reiten quiver, see for example \cite{MRS}. This correspondence was used to discover several connections between homological notions and combinatorial properties of Dyck paths, see for example \cite{MRS} and \cite{CM}. Motivated by this connection, we call a Dyck path $D$ \Dfn{shod} if the corresponding Nakayama algebra is shod and we call $D$ \Dfn{tilted} if the corresponding Nakayama algebra is tilted. In the preliminaries, we will discuss how to translate homological notions to elementary combinatorial notions on Dyck paths.
Our goal in this article is to give a combinatorial classification of shod Nakayama algebras using pattern avoiding permutations which are a classical topic in combinatorics with applications for example in algebraic geometry, see \cite{Bi} and \cite{Bo}.
Krattenthaler introduced in \cite{K} a bijection between Dyck paths and 132-avoiding permutations with many nice properties. We call this bijection the \Dfn{Krattenthaler bijection} in the following and refer to the preliminaries for the definition.
Our main result states the following:
\begin{theorem}
\label{thm::Krattenthaler_bij_intro}
The Krattenthaler bijection restricts to a bijection between shod Dyck paths and 132-avoiding permutations that additionally avoid the patterns 4321 and 4231.

This bijection further restricts to a bijection between tilted Dyck paths and 132-avoiding permutations that additionally avoid the pattern 321.
\end{theorem}
This combinatorial classification of shod Nakayama algebras via pattern avoiding permutations using the Krattenthaler bijection was first conjectured by John Machacek, \cite[Conjecture A]{Ma}.
As a corollary, we obtain a surprisingly beautiful formula for the enumeration of shod Nakayama algebras:
\begin{corollary}
\label{cor::counting_intro}
Let $n \geq 2$. 
\begin{enumerate}
\item The strict shod Nakayama algebras with $n$ simple modules are enumerated by $\sum\limits_{k=1}^{n-3}{k^2}$. 
\item The tilted Nakayama algebras with $n$ simple modules are enumerated by $1+\sum\limits_{k=1}^{n-2}{k}.$
\item The total number of shod Nakayama algebras with $n$ simple modules is given by $1+\sum\limits_{k=1}^{n-2}{k}+\sum\limits_{k=1}^{n-3}{k^2}=\frac{1}{3}(n^3-6n^2+14n-9)$.
\end{enumerate}
\end{corollary}

The next chapter introduces the required preliminaries on Nakayama algebras, their correspondance to Dyck paths, pattern avoidance and the Krattenthaler bijection. This allows us to prove our results in the following chapter. Finally, we give a brief outlook into possible generalisations.

\section{Preliminaries}

\subsection{Nakayama algebras}
In the following, $K$ always denotes a field and the algebras are finite-dimensional $K$-algebras. 
Recall that a finite-dimensional algebra $A$ is a Nakayama algebra if every indecomposable finite-dimensional $A$-module is uniserial, meaning that it has a unique composition series. 
We refer to the textbooks \cite{AF,ASS,SY,ARS} that all contain chapters on Nakayama algebras for a more detailed introduction to this important class of algebras and also the basics of representation theory of quiver algebras. 
If a Nakayama algebra $A$ is given by a quiver and relations, it is known that the underlying quiver is either $Q_n$, in which case we call it a \Dfn{linear Nakayama algebra}, or $\widetilde{Q}_n$ for $n \geq 1$ where $\widetilde{Q}_1$ is the one-loop quiver.

 \begin{equation*}
            Q_n: \begin{tikzcd}
0 \arrow[r, "a_1"] & 1 \arrow[r, "a_2"] & {...} \arrow[r, "a_{n-2}"] & n-2 \arrow[r, "a_{n-1}"] & n-1
\end{tikzcd}
\end{equation*}

\begin{equation*}
    \widetilde{Q}_n: \begin{tikzcd}
0 \arrow[r] & 1 \arrow[r] & {...} \arrow[r] & n-2 \arrow[r] & n-1 \arrow[llll, bend left]
\end{tikzcd}
\end{equation*}
    
Each Nakayama algebra is uniquely associated to a \Dfn{Kupisch series} $[c_0, c_1, \ldots, c_{n-1}]$ where $n$ coincides with the number of vertices of the underlying quiver and $c_i = \dim(P_i)$, the $K$-vector space dimension of the projective module $P_i = e_iA$.

\subsection{Dyck paths}

There is a well-known bijection which relates linear Nakayama algebras to Dyck paths which we will use throughout this paper. 
\begin{definition}
    Let us consider the grid $\mathbb{Z}^2$ written as in Figure \ref{fig::coord}. A \Dfn{Dyck path of semilength n} is a path from $(0,1)$ to $(n,1)$ consisting of up-steps $(0,1)$ and down-steps $(1,-1)$ such that the path never crosses the axis $y=1$. 
\end{definition}

To every Dyck path we can associate its so-called \Dfn{area sequence} $[c_0,c_1,\ldots,c_{n-1}]$ where $c_i$ is defined to be the maximal number such that $(i,c_i)$ is contained in the Dyck path. A Dyck path is uniquely determined by its area sequence. For example, the Dyck path in Figure \ref{fig::coord} has the area sequence $[3,2,1]$.

It is well-known that the defining relations of Dyck paths with semilength $n-1$ are the same as the defining relations of the Kupisch series of linear Nakayama algebras with $n$ simple modules. This implies the following bijection which has a central role in our paper, see \cite[Proposition 2.8]{MRS}.

\begin{proposition}
There is a bijection between linear Nakayama algebras with $n$ simple modules and Dyck paths of semilength $n-1.$ This is defined by sending a linear Nakayama algebra $A$ with Kupisch series $[c_0,c_1,\ldots,c_{n-1}]$ to the unique Dyck path $D_A$ of semilength $n-1$ which has area sequence equal to $[c_0,c_1,\ldots,c_{n-1}].$ 
\end{proposition}

\begin{remark}
    For algebraic purposes, we reshape the coordinate grid as in Figure \ref{fig::coord}. This way, the representation of the Dyck path corresponds to the Auslander-Reiten quiver of the respective Nakayama algebra.  Then every $\mathbb{Z}^2$-gridpoint enclosed by the Dyck path and the axis $y=1$ corresponds to an indecomposable $A$-module. From now on, we will always consider modules to be indecomposable and associate them with their respective coordinate point. 
\end{remark}

\begin{center}
\captionsetup{type=figure}
\includegraphics[]{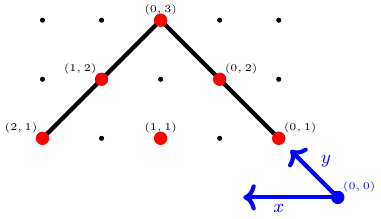}
\captionof{figure}{Coordinates of the indecomposable modules of the Nakayama algebra with Kupisch series $[3,2,1]$.}
\label{fig::coord}
\end{center}

The previous bijection is very useful as we can see many homological properties of the Nakayama algebra $A$ just by looking at the corresponding Dyck path $D_A.$ First, the indecomposable \Dfn{projective modules} correspond to the gridpoints on the left side of the Dyck path (the gridpoints marked green in Figure \ref{fig::proj_inj_mod}) and the indecomposable \Dfn{injective modules} correspond to the gridpoints on the right side of the Dyck path (the gridpoints marked blue in Figure \ref{fig::proj_inj_mod}). So the peaks of the Dyck path are exactly the modules which are both projective and injective. Note that valleys are neither injective nor projective.

\begin{center}
    \captionsetup{type=figure}
        \includegraphics[]{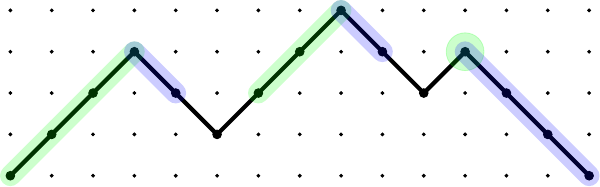}
        \captionof{figure}{Projective and injective modules.}
        \label{fig::proj_inj_mod}
    \end{center} 

Moreover, we can also describe the \Dfn{minimal projective resolution} of the modules in a completely combinatorial way on the Dyck path. For this, take any gridpoint $(i,j).$ Then we build its minimal projective resolution
$$\ldots \longrightarrow P_1 \longrightarrow P_0 \longrightarrow (i,j) \longrightarrow 0 $$
step-by-step.

In the first step, we obtain $P_0$ by starting at $(i,j)$ and taking steps in direction $(0,1)$ until we reach the left side of the Dyck path. So $P_0$ has coordinates $(i,c_i)$ where $[c_0,c_1,\ldots,c_{n-1}]$ is the area sequence of our Dyck path. Algebraically, this step corresponds to taking the \Dfn{projective cover} of the module $(i,j).$ If $(i,j)$ already lies on the left side of the Dyck path, i.e. is already projective then $P_0=(i,j).$ In the second step, we determine the gridpoint \textcolor{blue}{$\Omega(i,j)$}. We do this by starting at $P_0$ and taking $j$ steps in direction $(1,-1)$. This way we end up at $(i+j,c_i-j)=\Omega(i,j)$, see Figure \ref{fig::proj_inj_resol}. Algebraically, this corresponds to the \Dfn{(first) syzygy} module of $(i,j).$ 

Now, we will apply the first step to our new gridpoint $\Omega(i,j)$ to obtain $P_1.$ Then we apply the second step to calculate $\Omega^2(i,j)=\Omega(\Omega(i,j))$. By applying these two steps repeatedly, we obtain the entire projective resolution. This process stops if for some $k$ $\Omega^k(i,j)$ lies on the left side of the Dyck path, in other words, if  $\Omega^k(i,j)$ is projective. Then the minimal projective resolution has the form
$$0 \longrightarrow P_k = \Omega^k(i,j) \longrightarrow P_{k-1} \longrightarrow \ldots \longrightarrow P_1 \longrightarrow P_0 \longrightarrow (i,j) \longrightarrow 0. $$
The \Dfn{projective dimension} of $(i,j)$ \textcolor{blue}{$\pdim(i,j)$} is the minimal $k\ge 0$ such that the gridpoint $\Omega^k(i,j)$ is projective. 

\begin{center}
    \captionsetup{type=figure}
        \includegraphics[]{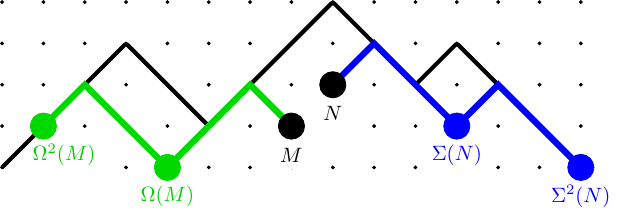}    
        \captionof{figure}{Projective and injective resolutions. Here $\pdim(M)=2$ and $\idim(N)=2$.}
        \label{fig::proj_inj_resol}
    \end{center} 
Dually, we can describe an algorithm to calculate a \Dfn{minimal injective resolution}, the \Dfn{injective envelope} (which is the dual notion to the projective cover), the \Dfn{(first) cosyzygy module} \textcolor{blue}{$\Sigma$}, and the \Dfn{injective dimension}, \textcolor{blue}{$\idim$} of a module corresponding to a gridpoint. 
This allows us to define the \Dfn{global dimension} of a Dyck path $D_A$ as the global dimension of its corresponding Nakayama algebra:
\begin{align*}
    \gldim D_A = \gldim A &= \sup \{ \pdim (i,j) \,|\, (i,j) \textnormal{ is an } A\textnormal{-module}\} \\
    &= \sup \{ \idim (i,j) \,|\, (i,j) \textnormal{ is an } A\textnormal{-module}\}.
\end{align*}

\subsection{Shod Nakayama algebras and Dyck paths}
On the level of Dyck paths, we are now able to understand all homological definitions needed to determine whether an algebra is shod. Recall that a finite-dimensional algebra is shod if $\pdim(M)\le 1$ or $\idim(M)\le 1$ for every module $M$. We say that a Dyck path D is \Dfn{shod} if the corresponding linear Nakayama algebra is shod. 

\begin{example}
    The Dyck path with area sequence $[4,3,3,2,1]$ is shod, see Figure \ref{fig::shod}.
    \begin{center}
    \captionsetup{type=figure}
        \includegraphics[]{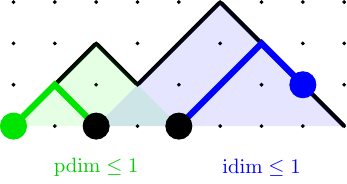}
        \captionof{figure}{Example of a shod Dyck path.}
        \label{fig::shod}
    \end{center} 
    The modules on the green side all have $\pdim \leq 1$, as symbolised by the module with the green projective resolution. Dually, the blue side contains modules with $\idim \leq 1$. 

    On the other hand, the Nakayama algebra with Kupisch series $[2,3,2,2,1]$ is not shod as there exists an indecomposable module $M$, marked orange in Figure \ref{fig::not_shod}, with projective and injective dimension 2.
    \begin{center}
    \captionsetup{type=figure}
        \includegraphics[]{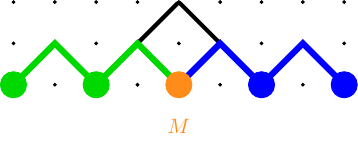}
        \captionof{figure}{Example of a non-shod Dyck path.}
        \label{fig::not_shod}
    \end{center}     
\end{example}

We will need the following results on shod algebras:
\begin{proposition}
\label{prop::tilted_fin_rep_type}
Let $A$ be an algebra of finite representation type.
Then $A$ is a quasi-tilted algebra if and only if $A$ is a tilted algebra.
\end{proposition} 
\begin{proof}
See \cite[Corollary 3.6]{HRS}.
\end{proof}
Since Nakayama algebras are always representation-finite, the notions quasi-tilted and tilted coincide for Nakayama algebras by the previous proposition.
\begin{proposition} \label{prop::nocycles}
A shod algebra has an acyclic quiver.

\end{proposition}
\begin{proof}
See \cite[Proposition 2.2]{CL}.

\end{proof}

In particular, for the classification of shod Nakayama algebras it is enough to restrict to linear Nakayama algebras by the previous proposition.

\subsection{Pattern avoiding permutations and the Krattenthaler bijection}

Lastly, we recall the \Dfn{Krattenthaler bijection $\phi$} introduced in \cite{K} which is a bijection from 132-avoiding permutations to Dyck paths. We will use this to describe a bijection between permutations avoiding certain patterns and shod Dyck paths. 

Let $\pi=\pi_1\pi_2\ldots\pi_n$ be a permutation of $\{1,\ldots,n\}$ and $\sigma=\sigma_1\sigma_2\ldots\sigma_m$ be a permutation of $\{1,\ldots,m\}$ for some $m\le n.$ We say that $\pi$ \Dfn{avoids the pattern} $\sigma$ if there are no indices $1\le i_1<i_2<\ldots<i_m\le n$ such that $\pi_{i_1}\pi_{i_2}\ldots\pi_{i_m}$ and $\sigma_1\sigma_2\ldots\sigma_m$ have the same relative order.

Now, let $\pi$ be a 132-avoiding permutation of length $n$. To describe $\phi(\pi)$, the Dyck path corresponding to $\pi$, we first have to determine all the left-to-right minima of $\pi$. A \Dfn{left-to-right minimum} of $\pi$ is an element $\pi_i$ such that it is smaller than all the elements to its left. In other words,  $\pi_i<\pi_j$ for every $j<i.$ Let us denote the left-to-right minima of $\pi$ by $m_1, m_2, \ldots, m_s.$ Then we can write $\pi$ as $$\pi=m_1w_1m_2w_2\ldots m_sw_s$$ where $w_i$ denotes the subword of $\pi$ between $m_i$ and $m_{i+1}.$ 

\begin{remark} 
\label{rmk::permutations}
Let $\pi=m_1w_1\ldots m_sw_s$ be a 132-avoiding permutation. The following are easy observations.
\begin{enumerate}[(i)]
    \item We have $m_s=1.$
    \item For every $1\le i\le  s$ $m_iw_i$ is strictly increasing. By the definition of the left-to-right minima, $m_i$ is smaller than every element of $w_i$. So if $w_i$ was not strictly increasing, $\pi$ would contain a $132$-pattern.
\end{enumerate}
\end{remark}

Now, $\phi(\pi)$ is defined to be the Dyck path which alternates between $|w_i|+1$ up-steps and $m_{i-1}-m_i$ down-steps.  Here $i$ ranges from $s$ to $1,$ in particular, the Dyck path starts with $|w_s|+1$ up-steps and finishes with $m_0-m_1$ down-steps. We set $m_0=n+1.$ Recall that, using our previous visualisation of grids, we read the Dyck paths from right to left. 

Note that $\phi(\pi)$ has semilength $(|w_1|+1)+(|w_2|+1)+\ldots +(|w_s|+1)=|\pi|=n.$ As a sanity check, we can also observe that $n=(m_0-m_1)+(m_1-m_2)+\ldots+(m_{s-1}-m_s).$ The action of this map $\phi$ is illustrated in Figure \ref{fig::Krattenthaler_bij}. 
\begin{center}
\captionsetup{type=figure}
    \includegraphics[]{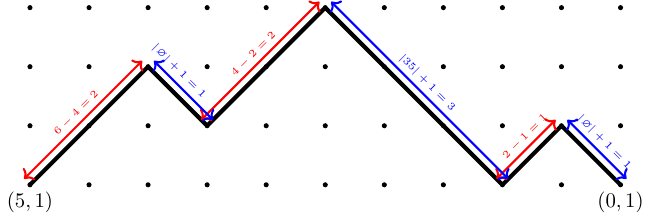}  
    \captionof{figure}{The Dyck path corresponding to $\pi = {\color{red}42}{\color{blue}35}{\color{red}1}$. We have $\color{red}m_1 = 4$, $\color{blue}w_1 = \varnothing$, $\color{red}m_2 = 2$, $\color{blue}w_2 = 35$, $\color{red}m_3 = 1$, $\color{blue}w_3 = \varnothing$ and $\pi = m_1w_1\ldots m_3w_3$.}
    \label{fig::Krattenthaler_bij}
\end{center}
Note that $s$, i.e. the number of left-to-right minima of $\pi$ is equal to the number of peaks of the Dyck path $\phi(\pi)$. 

\begin{proposition}\cite[Section 2]{K}
    The map $\phi$ is a well-defined bijection between 132-avoiding permutations of $\{1,2,\ldots,n\}$ and Dyck paths of semilength $n.$
\end{proposition}

\section{Shod Dyck paths}

\subsection{Classification}
This part aims to provide a classification of shod Dyck paths, both geometrically (see Theorem \ref{thm::strict_shod_d<1}) and in terms of pattern avoiding permutations (see Theorem \ref{thm::classificationKrattenthaler}). 

\begin{lemma}
\label{lem::atmost3peaks}
    Let $D$ be a shod Dyck path. Then it has at most three peaks.
\end{lemma}
\begin{proof}
    Let us assume for the sake of contradiction that $D$ has at least four peaks. Consider a module $M$ corresponding to a valley of $D$ enclosed by the peaks $Q_1$ and $Q_2$, such that $Q_1$ is not the leftmost peak and $Q_2$ is not the rightmost peak. This situation is visualised in Figure \ref{fig::4_peaks}.

    Assume that $\pdim(M)\le 1$. As $M$ is a valley, it is not projective, so the minimal projective resolution does not stop after the $0$-th term, and $\Omega(M)$ must be non-zero and projective. Following the rules described in the preliminaries, we see that the projective cover of $M$ is $Q_1$, and $\Omega(M)$ lies on the same line as $Q_1$ and $M',$ where $M'$ is the first valley to the left of $Q_1$.  $M'$ exists as $Q_1$ was chosen not to be the leftmost peak. Furthermore, $\Omega(M)$ is only projective if it lies strictly between $M'$ and $Q_1.$ This implies that the distance of $M'$ and $Q_1$, $d_1\ge 2.$ So the module with coordinates $(i+1,c_i-1)$ is projective. 
    
    Using this, we can show that the module $N$ with coordinates $(i+1,c_i-2)$ has projective cover $(i+1,c_i-1)$ and $\Omega(N)=(i+c_i-1,1),$ which is not projective as $Q_1$ is not the leftmost peak. This yields that $\pdim(N)\ge 2.$ Since $D$ was assumed to be shod, $\idim(N)\le 1$ has to hold. We know that $N$ is not injective, so $\Sigma(N)$ must be. If the distance of $M$ and $Q_1$, $d_2$ is at least $2$, then $\Sigma(N)=(i,1)$, which is not injective as $Q_1$ is not the rightmost peak. Thus, $d_2=1$ must hold. This, however, implies that $M$ has coordinates $(i,c_i-1)$, and so $\Omega(M) = \Omega(N)$. This is a contradiction to the assumption that $\pdim(M)\le 1.$ 

    \begin{center}
    \captionsetup{type=figure}
        \includegraphics[]{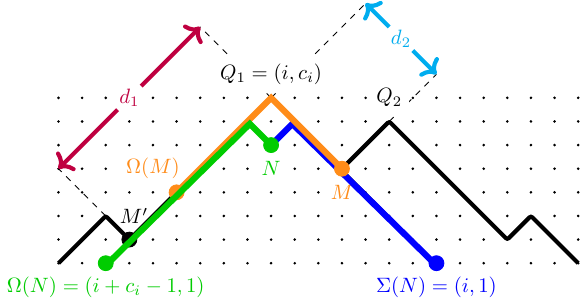}
        \captionof{figure}{Dyck paths with $\geq$ 4 peaks are not shod.}
        \label{fig::4_peaks}
    \end{center}

    A dual argument shows that $\idim(M)\le 1$ also leads to a contradiction. So if $D$ has at least four peaks then $\pdim(M)>1$ and $\idim(M)>1$ must hold. Therefore, $D$ cannot be shod.  
\end{proof}

The last lemma tells us that it is enough to concentrate on the Dyck paths with at most three peaks. In the next theorem we will use the following notation: We will denote the peaks from the left to the right by $Q_1, Q_2$ and $Q_3.$ We set $Q_2=Q_3$ if the Dyck path has only two peaks, and $Q_1=Q_2=Q_3$ if only one. We assume $Q_2$ to have coordinates $(i,c_i).$
We denote by $M_1$ the valley enclosed by $Q_1$ and $Q_2$, and by $M_2$ the valley enclosed by $Q_2$ and $Q_3.$ Furthermore, $d_i$ stands for the distance of $M_i$ and $Q_2$ for $i=1,2.$ This is illustrated in the following picture:

\begin{center}
\captionsetup{type=figure}
\includegraphics[]{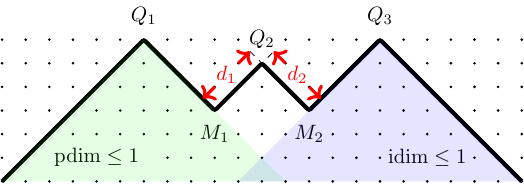}
    \captionof{figure}{Notation for three peaks.}
    \label{fig::notation_3_peaks}
\end{center}
The next theorem provides an elementary combinatorial classification of shod Nakayama algebras. We remark that this theorem can also be proven using the main result of \cite{BT} whose proof uses deep results from homological algebra and representation theory.
\begin{theorem}
\label{thm::strict_shod_d<1}
    Let $D$ be a Dyck path. Then it is shod if and only if it has at most two peaks or if it has three peaks and one of the following conditions holds:
    \begin{enumerate}[(A)]
        \item $d_1= 1$
        \item $d_2= 1.$
    \end{enumerate}
\end{theorem}
\begin{proof}
    Lemma \ref{lem::atmost3peaks} tells us that if $D$ has at least four peaks, $D$ is not shod. So it remains to classify which Dyck paths with at most three peaks are shod.
    
    It is easy to see that every gridpoint in the triangle with peak $Q_1$ (green triangle in Figure \ref{fig::notation_3_peaks}) has projective dimension at most $1,$ and every gridpoint in the triangle with peak $Q_3$ (blue triangle in Figure \ref{fig::notation_3_peaks}) has injective dimension at most $1.$ Moreover, every gridpoint lying on the Dyck path itself is either projective or injective. 
    
    If $D$ has at most two peaks or if it has three peaks and $d_1\le 1$ or $d_2\le 1$, then every gridpoint enclosed by $D$ is contained in one of these triangles or lies on the Dyck path, so $D$ must be shod. 

    On the other hand, if $D$ has three peaks but $d_1>1$ and $d_2>1$, then the gridpoint $N$ with coordinates $(i+1,c_i-2)$ (the orange gridpoint in the Figure \ref{fig::3_peaks_contradiction}) will violate the shod condition. 

    \begin{center}
    \captionsetup{type=figure}
\includegraphics[]{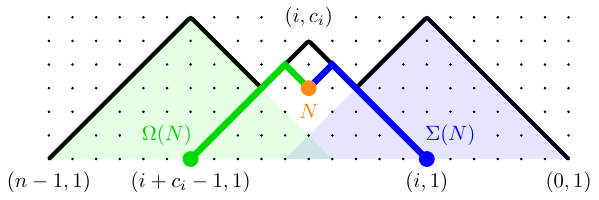}
    \captionof{figure}{The point $N$ violates the shod condition.}
    \label{fig::3_peaks_contradiction}
\end{center}
    
    Indeed, since $d_1>1,$ the module $(i+1,c_i-1)$ is the projective cover of $N,$ and $\Omega(N)=(i+c_i-1,1)$ (the green gridpoint in the picture), which is not projective. Thus, $\pdim(N)>1.$ Similarly, $d_2>1$ ensures that $(i,c_i-1)$ is the injective envelope of $N,$ and $\Sigma (N)= (i,1)$ (the blue gridpoint in the picture), which is not injective. Therefore, $\idim(N)>1.$ In conclusion, $D$ is not shod if it has at most three peaks and neither $(A)$ nor $(B)$ hold.  
\end{proof}

This classification restricts to a classification of tilted algebras.
\begin{lemma}
\label{lem::tilted}
A Dyck path is tilted if and only if it has at most two peaks.
\end{lemma}

\begin{proof}
    We already know by Proposition \ref{prop::tilted_fin_rep_type} that a Dyck path $D$ is tilted if and only if it is shod and $\gldim D \le 2$.
    Assume that $D$ has at most two peaks. It is easy to see that then the projective dimension of any gridpoint can be at most 2: Two projectives following directly after each other in the projective resolution are on the same section on the Dyck path if and only if the projective resolution stops afterwards. Since there are at most two sections of down-steps, the projective dimension is at most 2.

    Conversely, let $D$ be shod with three peaks and fulfilling condition $(A)$. We use the notation given by Figure \ref{fig::notation_3_peaks}. Let $Q_1 = (j,c_j)$, $Q_2 = (i,c_i)$. Observe that $(A)$ is equivalent to $j=i+1.$ The syzygy of the first valley $M_1 =(i+1, c_i-1)$ is given by $\Omega(M_1) = (i+c_i, c_j - c_i +1)$, see Figure \ref{fig::3_peaks_not_tilted}. Note that $\Sigma\Omega(M_1) = M_1$. Since $M_1$ is a valley, it is not injective. Moreover, $\Sigma^2\Omega(M_1)=\Sigma(M_1) = (i,1)$ which is not injective either as $Q_2$ is not the right-most peak. So the injective dimension of $\Omega(M_1)$ is at least 3. 
    \begin{center}
    \captionsetup{type=figure}
\includegraphics[]{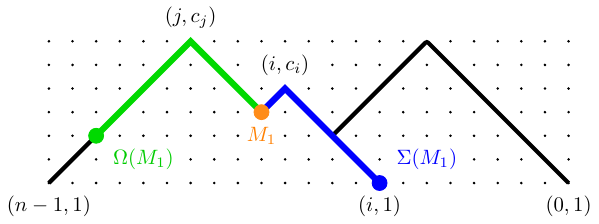}
    \captionof{figure}{The gridpoint $\Omega(M_1)$ has $\idim \ge 3$.}
    \label{fig::3_peaks_not_tilted}
\end{center}
    Then $3 \le \gldim D \le 3$. So the global dimension of $D$ is 3 and it is not tilted. Dually, if $(B)$ holds, we can work with the projective resolution of $\Sigma(M_2)$ to obtain the same result.
\end{proof}

Theorem \ref{thm::strict_shod_d<1} provides us with a characterisation of all shod Dyck paths. Interestingly, we can also characterise them in terms of the permutations they correspond to under the Krattenthaler bijection. This is what our main result describes.

\begin{theorem}
\label{thm::classificationKrattenthaler}
The Krattenthaler bijection restricts to a bijection between shod Dyck paths and permutations avoiding the patterns 132, 4321 and 4231. 
\end{theorem}

In order to prove this theorem, we first show the next lemma.  

\begin{lemma}
\label{lem::peaksAndPermutations}
    A Dyck path has at most $k$ peaks if and only if the corresponding $132$-avoiding permutation under the Krattenthaler bijection also avoids the pattern $(k+1)k\ldots1$. 
\end{lemma}
\begin{proof}
Assume $D$ is a Dyck path of semilength $n$, has at least $k+1$ peaks, and corresponds to the permutation $\pi=\pi_{1}\pi_{2}\ldots\pi_n=m_1w_1m_2w_2\ldots m_sw_s$. Note that $s$ is equal to the number of peaks, so $s\ge k+1.$ Moreover, for every $i$ $m_i>m_{i+1}$ since $m_{i+1}$ is the smallest element in $m_1w_1\ldots m_iw_im_{i+1}.$ Thus, $m_1>m_2>m_3>\ldots>m_{k+1}$ is a $(k+1)k\ldots1$-pattern in $\pi.$

Conversely, assume that $\pi$ contains a $(k+1)k\ldots1$ pattern; in other words, there are indices $1\le\sigma(1)<\sigma(2)<\ldots<\sigma (k+1)\le n$ such that $\pi_{\sigma(1)}>\pi_{\sigma(2)}>\ldots >\pi_{\sigma(k+1)}.$ Let $\pi_{\sigma(1)}$ be contained in $m_iw_i.$ 
Since $m_iw_i$ is strictly increasing by Remark \ref{rmk::permutations} (ii), $\pi_{\sigma(2)}$ must be contained in an $m_jw_j$ with $j>i.$ Repeating this argument, we obtain an increasing sequence $i_1<i_2<\ldots<i_{k+1}$ such that $\pi_{i_{\ell}}\in m_{i_{\ell}}w_{i_{\ell}}$ for every $1\le \ell \le k+1.$ This implies that there are at least $k+1$ $m_i$-s, hence, $k+1$ peaks.  

\end{proof}

\begin{proof}[Proof of Theorem 3.3.]
Let $D$ be a Dyck path with at most three peaks of semilength $n$ and $\pi$ the permutation corresponding to $D$ under the Krattenthaler bijection. By Lemma \ref{lem::atmost3peaks} and \ref{lem::peaksAndPermutations}, it remains to show that $D$ is shod if and only if $\pi$ avoids the pattern $4231$. 

Assume $D$ is shod. Then by Theorem \ref{thm::strict_shod_d<1}, it has either at most two peaks, or it has three peaks and condition $(A)$ or $(B)$ from the same theorem is fulfilled. By Lemma \ref{lem::peaksAndPermutations}, if $D$ has at most two peaks, it avoids the pattern $321,$ and in particular, it also avoids the pattern $4231$. 

Since the number of left-to-right minima of $\pi$ is equal to the number of peaks of the corresponding Dyck path, if $D$ has three peaks, then $\pi$ has the form $$\pi=m_1w_1m_2w_2m_3w_3.$$ Figure \ref{fig::Krattenthaler_bij_3_peaks} illustrates $D.$ 
\begin{center}
\captionsetup{type=figure}
    \includegraphics[]{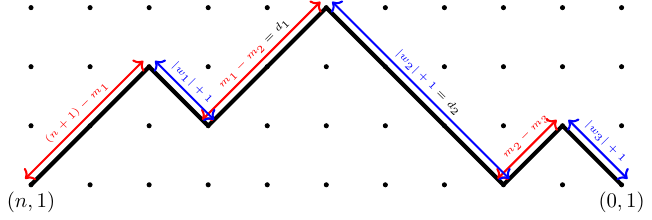}
    \captionof{figure}{The Dyck path corresponding to $\pi = {\color{red}m_1}{\color{blue}w_1}{\color{red}m_2}{\color{blue}w_2}{\color{red}m_3}{\color{blue}w_3}$.}
    \label{fig::Krattenthaler_bij_3_peaks}
\end{center}

Assume for the sake of contradiction that there are indices $\sigma(1)<\sigma(2)<\sigma(3)<\sigma(4)$ such that $\pi_{\sigma(4)}<\pi_{\sigma(2)}<\pi_{\sigma(3)}<\pi_{\sigma(1)}$. As $m_iw_i$ is increasing by Remark \ref{rmk::permutations} (ii), $\pi_{\sigma(1)} \in m_1w_1$, $\pi_{\sigma(2)},\pi_{\sigma(3)} \in m_2w_2$ and $\pi_{\sigma(4)} \in m_3w_3$. Since $\sigma(2)<\sigma(3)$, $\pi_{\sigma(3)}\in w_2.$

First, we look at the case when $D$ satisfies $(B).$ Observe that $d_2=1$ is equivalent to $|w_2|=0.$  This contradicts to $\pi_{\sigma(3)}\in w_2.$ 

Second, assume that $D$ satisfies $(A)$ but does not satisfy $(B).$ Note that this is equivalent to $m_1=m_2+1$ and $|w_2|\ge 1.$ Since every element of $w_2$ is bigger than $m_2,$ it is also bigger than $m_1=m_2+1,$ so  $\pi_{\sigma(3)}>m_1$. Then $\pi_{\sigma(1)}>\pi_{\sigma(3)}$ yields $\pi_{\sigma(1)}\in w_1.$

Observe that every element of $w_1$ must be smaller than every element in $w_2,$ otherwise there would be a $132$-pattern in $\pi.$ Indeed, choose any $v\in w_2.$ We already saw that $v>m_1.$ If there is an element $u\in w_1$ with $u>v,$ then $m_1uv$ gives a $132$-pattern in $\pi,$ which is impossible. This, however, means that there is no $v=\pi_{\sigma(3)}\in w_2$ such that $\pi_{\sigma(3)}<\pi_{\sigma(1)}=u$. So we reached a contradiction in both cases allowing us to conclude that if $D$ is shod, $\pi$ avoids the pattern $4231.$

For the other direction, assume that $D$ is not shod. Then, by Theorem \ref{thm::strict_shod_d<1} it satisfies neither $(A)$ nor $(B).$ We saw that this translates to $m_1>m_2+1$ and $|w_2|\neq 0.$ We will show that $m_2+1\in w_2$, which yields that $m_1m_2(m_2+1)m_3$ forms a $4231$-pattern in $\pi.$  Since $m_1w_1$ is strictly increasing, $m_2+1<m_1$ cannot be in $m_1w_1.$ We assumed that $|w_2|\neq 0.$ If $m_2+1\notin w_2$, there exists a $v\in w_2$ such that $v>m_2+1$ and $m_2+1\in w_3.$ But then, $m_2v(m_2+1)$ would correspond to a $132$-pattern in $\pi$ which is impossible. Hence, $m_2+1$ is necessarily in $w_2.$ This concludes the proof of the theorem.    
\end{proof}

\begin{corollary}
    The Krattenthaler bijection restricts to a bijection between tilted Dyck paths and permutations avoiding the patterns $132$ and $321$.
\end{corollary}
\begin{proof}
This follows immediately from Lemma \ref{lem::tilted} and \ref{lem::peaksAndPermutations}. As a sanity check, if a pattern avoids $321$ then it must be avoiding $4321$ and $4231$ as well. So this is in line with the statement of Theorem \ref{thm::classificationKrattenthaler}.
\end{proof}

\subsection{Enumeration}

With the complete classification of shod Dyck paths, and thus, shod Nakayama algebras, we are now able to study how they are enumerated.
    
\begin{lemma}
\label{lem::counting_tilted}
    There are $1 + \sum_{k = 1}^{n-2} k$ tilted Nakayama algebras with $n$ simples. Equivalently, this enumerates the number of tilted Dyck paths of semilength $n-1$, or permutations of $\{1, \ldots, n-1\}$ avoiding the patterns $132$ and $321$.
\end{lemma}

\begin{proof}
    The tilted Nakayama algebras correspond to Dyck paths of semilength $n-1$ with one or two peaks. There is only one Dyck path with one peak. Let us consider a Dyck path with two peaks with $k$ denoting the distance of $(n-1,1)$ to the first peak and $l$ the one of $(0,1)$ to the last peak, see Figure \ref{fig::counting_tilted}.
    \begin{center}
    \captionsetup{type=figure}
\includegraphics[]{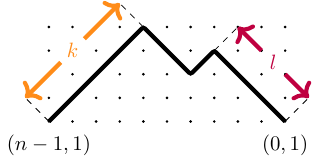}
\captionof{figure}{Counting two peaks.}
\label{fig::counting_tilted}
\end{center}
    We know that $1 \le k \le n-2$ since the $n-1$ down-steps are split into two non-empty intervals. Then the other interval of down-steps has length $n-1-k$. Thus, $l$ must be at least $n-1-k$ in order for the Dyck path to stay above the axis $y = 1$. Moreover, $l$ is also bounded above by $n-2$. So the number of Dyck paths with two peaks is 
\begin{equation*}
    \sum_{k=1}^{n-2} \sum_{l=n-1-k}^{n-2} 1 \, = \sum_{k = 1}^{n-2} k.
\end{equation*}

\end{proof}

\begin{lemma}
\label{lem::counting_strict}
    There are $\sum_{k = 1}^{n-3} k^2$ strict shod linear Nakayama algebras with $n$ simples. Equivalently, this enumerates the number of strict shod Dyck paths of semilength $n-1$.
\end{lemma}

\begin{proof}
Using Theorem \ref{thm::strict_shod_d<1}, it is enough to count Dyck paths with three peaks where $d_1 = 1$ or $d_2 = 1$. Let us consider the case $d_1 = 1$, the other is counted dually. Let $k$ be the distance of $(n-1,1)$ to the first peak and $l$ the one of $(0,1)$ to the last peak, see Figure \ref{fig::counting}.
\begin{center}
\captionsetup{type=figure}
\includegraphics[]{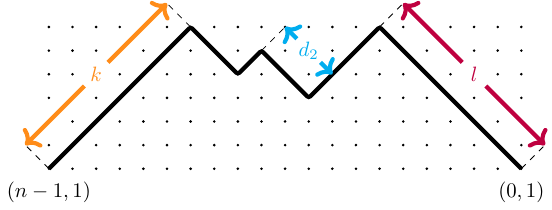}
\captionof{figure}{Counting in the case $d_1 = 1$.}
\label{fig::counting}
\end{center}
We know that $1 \leq k \leq n-3$ since the $n-1$ down-steps have to be split into three non-empty sections. As $d_1 = 1$, this implies that the distance from the second valley to the third peak is $n-2-k$. So $l$ must be at least $n-2-k$, otherwise the Dyck path would go below the axis $y=1$. Like $k$, $l$ too must not be larger than $n-3$. Finally, $d_2$ may range between 1 and $n-l-2$.

Since the case of $d_2 = 1$ is a mirrored version of the case $d_1 = 1$, they both account for the same number of strict shod Dyck paths. So to get the total number of strict shod Dyck paths, we only need to consider which cases we count twice. These are the ones where $d_1 = 1 = d_2$. In conclusion, the number of Dyck paths with three peaks where $d_1 = 1$ or $d_2 = 1$ is
\begin{align*}
    \sum_{k=1}^{n-3} \sum_{l=n-2-k}^{n-3} \left( 2 \left(\sum_{d_2=1}^{n-l-2} 1 \right) - 1 \right)
    &= \sum_{k=1}^{n-3} \sum_{l=n-2-k}^{n-3} \left( 2 \left( n-2-l \right) - 1 \right) \\
    &= \sum_{k=1}^{n-3} \sum_{l'=1}^{k} \left( 2 l' - 1 \right) = \sum_{k=1}^{n-3} k^2.   
\end{align*}
\end{proof}

Thus, Lemma \ref{lem::counting_tilted} and \ref{lem::counting_strict} combined prove Corollary \ref{cor::counting_intro}.

\section{Outlook}
We expect that certain higher dimensional generalisations of shod Nakayama algebras can also be described using pattern avoiding permutations.
Shod Nakayama algebras correspond to Dyck paths that avoid a certain homological pattern, namely indecomposable modules with projective dimension at least 2 and injective dimension at least 2.
We give some definitions generalising shod algebras that can be seen as some sort of pattern avoidance for Dyck paths when specialised to Nakayama algebras. 
We give the following generalisation for general finite-dimensional algebras:
\begin{definition}
Let $A$ be a finite-dimensional algebra over a field $K$ and $n,m$ two natrual numbers. Then $A$ is called \Dfn{$(n,m)$-shod} if for every indecomposable $A$-module $M$ we have $\pdim M \leq n$ or $\idim M \leq m$.
$A$ is called \Dfn{strong $(n,m)$-shod} if it is $(n,m)$-shod and additionally $A$ has global dimension at most $n+m$. \end{definition}
Note that for $n=m=1$, (1,1)-shod algebras are exactly the classical shod algebras and strong (1,1)-shod algebras are exactly the quasi-tilted algebras. As in \cite[Proposition 1.1, Chapter II]{HRS} one can show that the global dimension of an $(n,m)$-shod algebra is bounded above by $n+m+1$.
We pose the following problem:
\begin{problem}
Classify the (n,m)-shod and strong (n,m)-shod linear Nakayama algebras.
\end{problem}

We obtained some conjectures for small values of $n$ and $m$ that we plan to examine in forthcoming work. The problem is, however, much more complicated than the one for the classical shod algebras since the Krattenthaler bijection does usually not give us directly a correspondance between $(n,m)$-shod Dyck paths and a nice class of pattern avoiding permutations.

\printbibliography

\end{document}